\documentclass [12pt] {article}
\def\to{\rightarrow}
\begin{document}
\title{An asymptotic expansion for a ratio of products of gamma functions}
\author{Wolfgang B\"uhring\\
{\it Physikalisches Institut, Universit\"at Heidelberg, Philosophenweg 12,}\\ 
{\it 69120 Heidelberg, GERMANY}}

\date{}
\maketitle
\begin{abstract}
An asymptotic expansion of a ratio of products of gamma functions is derived.  It generalizes a formula which was stated by Dingle, first proved by Paris, and recently reconsidered by Olver.
\end{abstract}



\bigskip
 buehring@physi.uni-heidelberg.de

\section{Introduction}

Our starting point is the Gaussian hypergeometric function $F(a,b;c;z)$ and its series representation
\[
{1 \over {\Gamma (c)}}F(a,b;c;z)=\sum\limits_{n=0}^\infty  {{{(a)_n(b)_n} \over {\Gamma (c+n){n! }}}{{z^n} }},\quad |z|<1,
\]
which here is written in terms of Pochhammer symbols
\[(x)_n=x(x+1)\ldots \,(x+n-1)=\Gamma (x+n)/\Gamma (x).\]
The hypergeometric series appears as one solution of the Gaussian (or hypergeometric) differential equation, which is characterized by its three regular singular points at $z=0, 1, \infty$. The local series solutions at $0$ and $1$ of this differential equation are connected by the continuation formula \cite{abra}
\[{1 \over {\Gamma (c)}}F(a,b;c;z)={{\Gamma (c-a-b)} \over {\Gamma (c-a)\Gamma (c-b)}}F(a,b;1+a+b-c;1-z)\]
\begin{equation}
\mbox{}+{{\Gamma (a+b-c)} \over {\Gamma (a)\Gamma (b)}}(1-z)^{c-a-b}F(c-a,c-b;1+c-a-b;1-z),
\label{e2}
\end{equation}
\[(|\arg (1-z)|<\pi) .\]

Here we want to show that Eq.\ (\ref{e2}) implies an interesting asymptotic expansion for a ratio of products of gamma functions, of which only a special case was known before.

By applying the method of Darboux \cite{olv0,ss} to (\ref{e2}), we derive in Sec.\ \ref{sec:2} the formula in question. The behaviour of this and a related formula is discussed in Sec.\ \ref{sec:3} and illustrated by a few numerical examples.

\section{Derivation of an asymptotic expansion for a ratio of products of gamma functions}
\label{sec:2}

It is well-known that the late coefficients of a Taylor series expansion contain information about the nearest singular point of the expanded function \cite{mofe}. In this respect we want to analyze the continuation formula (\ref{e2}), in which then only the second, at $z=1$ singular term $R$ is relevant, which may be written as
\[R={{\Gamma (a+b-c)\Gamma (1+c-a-b)} \over {\Gamma (a)\Gamma (b)}}\sum\limits_{m=0}^\infty  {{{(c-a)_m(c-b)_m} \over {\Gamma (1+c-a-b+m)m! }}(1-z)^{c-a-b+m}}.\]
By means of the binomial theorem in its hypergeometric-series-form , we may expand the power factor
\[(1-z)^{c-a-b+m}=\sum\limits_{n=0}^\infty  {{{\Gamma (a+b-c-m+n)} \over {\Gamma (a+b-c-m)n! }}}z^n.\]
Interchanging then the order of the summations and simplifying by means of the reflection formula of the gamma function, we arrive at 
\[R={1 \over {\Gamma (a)\Gamma (b)}}\sum\limits_{n=0}^\infty  {\sum\limits_{m=0}^\infty  {(-1)^m{{(c-a)_m(c-b)_m} \over {m! }}{{\Gamma (a+b-c-m+n)} \over {n! }}z^n}}.\]
This is to be compared with the left-hand side $L$ of (\ref{e2}), which is
\[L={1 \over {\Gamma (a)\Gamma (b)}}\sum\limits_{n=0}^\infty  {{{\Gamma (a+n)\Gamma (b+n)} \over {\Gamma (c+n)n! }}z^n}.\]
Comparison of the coefficients of these two power series, which according to Darboux \cite{olv0} and Sch\"afke and Schmidt \cite{ss} should agree asymptotically as $n\to\infty$, then yields
\begin{equation}
\label{e3}
{{\Gamma (a+n)\Gamma (b+n)} \over {\Gamma (c+n)}}=\sum\limits_{m=0}^M {(-1)^m{{(c-a)_m(c-b)_m} \over {m! }}\Gamma (a+b-c-m+n)}
\end{equation}
\[\mbox{}+O(\Gamma (a+b-c-M-1+n)).\]
By means of
\[O(\Gamma (a+b-c-M-1+n))=\Gamma (a+b-c+n)O(n^{-M-1})\]
and the reflection formula of the gamma function, the relevant formula (\ref{e3}) may also be written as
\begin{equation}
\label{e4}
{{\Gamma (a+n)\Gamma (b+n)} \over {\Gamma (c+n)\Gamma (a+b-c+n)}}=1+\sum\limits_{m=1}^M {{{(c-a)_m(c-b)_m} \over {m! (1+c-a-b-n)_m}}}+O(n^{-M-1}).
\end{equation}
The asymptotic expansion for a ratio of products of gamma functions in this form (\ref{e4}) or the other (\ref{e3}) seems to be new. It is only the special case when  $c = 1$ which is known. This special case was stated by Dingle\cite{din}, first proved by Paris\cite{par},  and reconsidered recently by Olver\cite{olv1}, who has found a simple direct proof. His proof, as well as the proof of Paris, can be adapted easily to the more general case when   $c$   is different from   $1$ . Still another proof is available \cite{olv2} which includes an integral representation of the remainder term.
Our derivation of Eq. (\ref{e3}) or (\ref{e4}) is significantly different from all the earlier proofs of the case when $c=1$.

\section{Discussion and numerical examples}
\label{sec:3}

We now want to discuss our result in the form (\ref{e4}). First we observe that  the substitution $c\to a+b-c$  leads to the related formula
\begin{equation}
\label{e5}
{{\Gamma (a+n)\Gamma (b+n)} \over {\Gamma (c+n)\Gamma (a+b-c+n)}}=1+\sum\limits_{m=1}^M {{{(a-c)_m(b-c)_m} \over {m! (1-c-n)_m}}}+O(n^{-M-1}).
\end{equation}
Which of (\ref{e4}) or (\ref{e5}) is more advantageous numerically depends on the values of the parameters, and in this respect the two formulas complement each other. Table \ref{t1} shows an example with a set of parameters for which (\ref{e4}) gives more accurate values than (\ref{e5}), while Table \ref{t2} contains an example for which (\ref{e5}) is superior to (\ref{e4}). 

For  finite $n$ and $M \to\infty$  the series on the right-hand side of (\ref{e4}) converges if  $\hbox{Re}(1-c-n)>0$. The same is true for (\ref{e5}) if  $\hbox{Re}(1+c-a-b-n)>0$ . Then, in both cases, the Gaussian summation formula yields
\[{{\Gamma (1-c-n)\Gamma (1+c-a-b-n)} \over {\Gamma 1-a-n)\Gamma (1-b-n)}},\]
which, by means of the reflection formula of the gamma function, is seen to be equal to
\begin{equation}
\label{e6}
{{\Gamma (a+n)\Gamma (b+n)} \over {\Gamma (c+n)\Gamma (a+b-c+n)}}{{\sin (\pi [a+n])\sin (\pi [b+n])} \over {\sin (\pi [c+n])\sin (\pi [a+b-c+n])}}.
\end{equation}
Otherwise (\ref{e3}) -- (\ref{e5}) are divergent asymptotic expansions as $n\to\infty$. 

Although in our derivation  $n$   is a sufficiently large positive integer, the asymptotic expansions (\ref{e3}) -- (\ref{e5}) are expected to be valid in a certain sector of the complex  $n$  -plane, and in fact, the proofs of Paris \cite{par} and of Olver \cite{olv2} apply to complex values of  $n$. 

If the series in (\ref{e4}) or (\ref{e5}) converge, their sums are equal to (\ref{e6}), which generally (if neither $c-a$ nor $c-b$ is equal to an integer ) is different from the left-hand side of (\ref{e4}) or (\ref{e5}). Therefore (\ref{e4}) and (\ref{e5}) can be valid only in the half-planes in which the series do not converge. This means that (\ref{e4}) is an asymptotic expansion as  $n\to\infty$  in the half-plane $\hbox{Re}(c-1+n) \geq 0$, and (\ref{e5}) is an asymptotic expansion as  $n\to\infty$  in the half-plane $\hbox{Re}(a+b-c-1+n) \geq 0$. Otherwise the series on the right-hand sides represent a different function, namely (\ref{e6}).
 
A few numerical examples may serve for demonstration of these facts. In Table \ref{t3} , the series converge to (\ref{e6}) for $n = 10$  , and therefore (\ref{e4}) and (\ref{e5}) are not valid. For $ n = 20$, on the other hand, the series diverge and so (\ref{e4}) and (\ref{e5})  hold. The transition between the two regions is at the line $\hbox{Re}(n )= 12.4 $ in case of (\ref{e4}) or  $\hbox{Re}(n )= 12.5$  in case of (\ref{e5}). In Table \ref{t4}, we see convergence for $n = -15$   and divergence for  $n = -5$, the transition between the two regions being at the line  $\hbox{Re}(n )= -10.4$  in case of (\ref{e4}) or  $\hbox{Re}(n )= -10.5$  in case of (\ref{e5}) .

\newpage
\begin{table}
\begin{tabular}{ccrlrl}
 &$M$  &\multicolumn{2}{l}{right-hand 
side of (\ref{e4})}&\multicolumn{2}{l}{right-hand 
side of(\ref{e5})}\\ \hline
$n=10$  & 	 1	&0.9771429&&	0.9744681&  \\
	& 2&	0.9773113& &	0.9780243& \\
	& 3&	0.9772978&&	0.9769927&\\
	 &4&	0.9773005&&	0.9774980 & $\leftarrow$\\
	 &5&	0.9772995 & $\leftarrow$&	0.9771117 &$\leftarrow$\\
	 &6	&0.9773001& $\leftarrow$	&0.9775615& \\
	 &7&	0.9772995&&	0.9767519& \\
	&8&	0.9773003&&	0.9791530&\\
	& 9&	0.9772983&&	0.9652341&\\
	&10&	0.9773079&&	1.2823765&\\ 
&\multicolumn{5}{l}{exact value of  (\ref{e4}) or (\ref{e5}):  0.97729983 }\\
\end{tabular}
\caption{Values of the right-hand sides of (\protect\ref{e4}) and (\protect\ref{e5}) for the parameters
	$a = 0.7$,      $ b = 1.2$,       $c =  0.4$. }
\label{t1}
\end{table}

\begin{table}
\begin{tabular}{ccrlrl}
 & $M$  &\multicolumn{2}{l}{right-hand 
side of (\ref{e4})}&\multicolumn{2}{l}{right-hand 
side of(\ref{e5})}\\ \hline

$n=10$  &  1&	0.968000&&	0.972093&\\
	& 2&	0.973760&&	0.972350&\\
	& 3&	0.971512&$\leftarrow$ &	0.972324 &\\
	& 4&	0.973078 & $\leftarrow$ &	0.972331 &\\
	& 5 &	0.971231&&	0.972327 & $\leftarrow$ \\
	 &6&	0.975016&&	0.972330 &$\leftarrow$ \\
	& 7&	0.959571&&	0.972325 &\\
	& 8	&1.179434&&	0.972342 &\\
	& 9	&4.748048&&	0.972163&\\
	&10  &26.430946&&	0.968966&\\
&\multicolumn{5}{l}{exact value  of (\ref{e4}) or (\ref{e5}):   0.97232850}\\	
\end{tabular}
\caption{Values of the right-hand sides of (\protect\ref{e4}) and (\protect\ref{e5}) for the parameters
$a = -0.7$,        $b = -1.2$,       $c =  -0.4$.  }
\label{t2}
\end{table}

\begin{table}
\begin{tabular}{ccrlrl}
 & $M$  &\multicolumn{2}{l}{right-hand 
side of (\ref{e4})}&\multicolumn{2}{l}{right-hand 
side of(\ref{e5})}\\ \hline
$n=10$ &  	 1&	0.976000&&	0.975000&\\
& 2&	0.972434&&	0.971912&\\
&	 3&	0.971341&&	0.971037&\\
	& 4	&0.970882&&	0.970687&\\
& 5	&0.970651&&	0.970517&\\
&	 6&	0.970520&&	0.970423&\\
&	 7&	0.970440&&	0.970367&\\
	& 8&	0.970388&&	0.970331&\\
&	 9&	0.970352&&	0.970307&\\
&	10	&0.970326&&	0.970290&\\
&\multicolumn{5}{l}{exact value  of(\ref{e4}) or (\ref{e5}):  1.94045281}\\
&\multicolumn{5}{l}{exact value of (\ref{e6}):  0.97022640  \quad $\leftarrow$}\\ \hline

$n=20$ &  	 1&	1.008000&&	1.007895&\\
	&2&	1.007360&&	1.007392&\\
	& 3&	1.007521&&	1.007504&\\
	& 4 &	1.007438  &$\leftarrow$& 1.007452   &$\leftarrow$\\
	& 5&	1.007515  &$\leftarrow$  &	1.007497  &$\leftarrow$\\
	& 6&	1.007385&&	1.007426&\\
	& 7&	1.007839&&	1.007650&\\
	& 8	&1.002201&&	1.005398&\\
	 &9&	0.921096&&	0.965891&\\
	&10	&0.478588&&	0.740024&\\
&\multicolumn{5}{l}{	exact value  of (\ref{e4}) or (\ref{e5}):  1.00747290 \quad  $\leftarrow$}\\
&\multicolumn{5}{l}{exact value of (\ref{e6}):  0.50373645}\\ 
\end{tabular}
\caption{Values of the right-hand sides of (\protect\ref{e4}) and (\protect\ref{e5}) for the parameters
$a =-11.7$,        $ b = -11.2$,        $c = -11.4$ . } 
\label{t3}
\end{table}
  
\begin{table}
\begin{tabular}{ccrlrl}
 & $M$  &\multicolumn{2}{l}{right-hand 
side of (\ref{e4})}&\multicolumn{2}{l}{right-hand 
side of(\ref{e5})}\\ \hline
$n=-15$ &  	 1&	0.986667&&	0.986957&\\	
	& 2&	0.985648	&&0.985745&\\
	& 3&	0.985453&&	0.985492&\\
	& 4&	0.985397&&	0.985415&\\
	& 5&	0.985376&&	0.985386&\\
	& 6&	0.985368&&	0.985373&\\
	&7&	0.985363&&	0.985367&\\
	& 8&	0.985361&&	0.985363&\\
	& 9&	0.985360&&	0.985361&\\
	&10	&0.985359&&	0.985360&\\
	&\multicolumn{5}{l}{exact value of (\ref{e4}) or (\ref{e5}):  1.97071532}\\
	&\multicolumn{5}{l}{exact value of (\ref{e6}):  0.98535766 \quad $\leftarrow$}\\  \hline

$n=-5$ &	 1	&1.010909&&1.011111&\\
	 &2	&1.009891&&	1.009798&\\
	& 3&	1.010254 & $\leftarrow$ &1.010331  & $\leftarrow$\\
	& 4	&1.009940& $\leftarrow$&1.009818 & $\leftarrow$\\
	& 5	&1.010589	&&1.011015 &\\ 
	& 6&	1.005300&&	0.998322&\\
	& 7&	0.951894&&0.887892&\\
	 &8	&0.737202 &&	0.459630 &\\                         
	& 9&	0.134729 && $-$0.725230&\\
	&10 &$-$1.243041 &&  $-$3.418810&\\
	&\multicolumn{5}{l}{exact value  of (\ref{e4}) or (\ref{e5}):  1.01011438  \quad $\leftarrow$}\\
	&\multicolumn{5}{l}{exact value of (\ref{e6}):  0.50505719}\\ 
\end{tabular}
\caption{Values of the right-hand sides of (\protect\ref{e4}) or (\protect\ref{e5}) for the parameters
	$a = 11.7$,       $ b = 11.2$,      $c = 11.4$. } 
\label{t4}    
\end{table}
\end{document}